# On Calderón's conjecture

By MICHAEL LACEY and CHRISTOPH THIELE*

## 1. Introduction

This paper is a successor of [4]. In that paper we considered bilinear operators of the form

$$H_\alpha(f_1, f_2)(x) := \text{p.v.} \int f_1(x-t) f_2(x + \alpha t) \frac{dt}{t}, \tag{1}$$

which are originally defined for $f_1$, $f_2$ in the Schwartz class $\mathcal{S}(\mathbb{R})$. The natural question is whether estimates of the form

$$\|H_\alpha(f_1, f_2)\|_p \le C_{\alpha, p_1, p_2} \|f_1\|_{p_1} \|f_2\|_{p_2} \tag{2}$$

with constants $C_{\alpha,p_1,p_2}$ depending only on $\alpha, p_1, p_2$ and $p := \frac{p_1 p_2}{p_1 + p_2}$ hold. The first result of this type is proved in [4], and the purpose of the current paper is to extend the range of exponents $p_1$ and $p_2$ for which (2) is known. In particular, the case $p_1 = 2$, $p_2 = \infty$ is solved to the affirmative. This was originally considered to be the most natural case and is known as Calderón's conjecture [3].

We prove the following theorem:

THEOREM 1. *Let $\alpha \in \mathbb{R} \setminus \{0, -1\}$ and*

$$1 < p_1, p_2 \le \infty, \tag{3}$$

$$\frac{2}{3} < p := \frac{p_1 p_2}{p_1 + p_2} < \infty. \tag{4}$$

*Then there is a constant $C_{\alpha,p_1,p_2}$ such that estimate* (2) *holds for all $f_1, f_2 \in \mathcal{S}(\mathbb{R})$.*

If $\alpha = 0, -1, \infty$, then we obtain the bilinear operators

$$H(f_1) \cdot f_2, \; H(f_1 \cdot f_2), \; f_1 \cdot H(f_2),$$

*The first author acknowledges support from the NSF. Both authors have been supported by NATO travel grants and the Volkswagen Stiftung (RiP program in Oberwolfach).

the last one by replacing $t$ with $t/\alpha$ and taking a weak limit as $\alpha$ tends to infinity. Here $H$ is the ordinary linear Hilbert transform, and $\cdot$ is pointwise multiplication. The $L^p$-bounds of these operators are easy to determine and quite different from those in the theorem. This suggests that the behaviour of the constant $C_{\alpha,p_1,p_2}$ is subtle near the exceptional values of $\alpha$. It would be of interest to know that the constant is independent of $\alpha$ for some choices of $p_1$ and $p_2$.

We do not know that the condition $\frac{2}{3} < p$ is necessary in the theorem. But it is necessary for our proof. An easy counterexample shows that the unconditionality in inequality (6) already requires $\frac{2}{3} \leq p$. The cases of $(p_1, p_2)$ being equal to $(1, \infty)$, $(\infty, 1)$, or $(\infty, \infty)$ have to be excluded from the theorem, since the ordinary Hilbert transform is not bounded on $L^1$ or $L^\infty$.

We assume the reader as somewhat familiar with the results and techniques of [4]. The differences between the current paper and [4] manifest themselves in the overall organization and the extension of the counting function estimates to functions in $L^q$ with $q < 2$.

The authors would like to thank the referee for various corrections and suggestions towards improving this exposition.

## 2. Preliminary remarks on the exponents

Call a pair $(p_1, p_2)$ good, if for all $\alpha \in \mathbb{R} \setminus \{0, -1\}$ there is a constant $C_{\alpha,p_1,p_2}$ such that estimate (2) holds for all $f_1, f_2 \in \mathcal{S}(\mathbb{R})$. In this section we discuss interpolation and duality arguments. These, together with the known results from [4], show that instead of Theorem 1 it suffices to prove:

PROPOSITION 1. *If $1 < p_1, p_2 < 2$ and $\frac{2}{3} < \frac{p_1 p_2}{p_1 + p_2}$, then $(p_1, p_2)$ is good.*

In [4] the following is proved:

PROPOSITION 2. *If $2 < p_1, p_2 < \infty$ and $1 < \frac{p_1 p_2}{p_1 + p_2} < 2$, then $(p_1, p_2)$ is good.*

Strictly speaking, this proposition is proved in [4] only in the case $\alpha = 1$, but this restriction is inessential. The necessary modifications to obtain the full result appear in the current paper in Section 3. Therefore we take Proposition 2 for granted.

The next lemma follows by complex interpolation as in [1]. The authors are grateful to E. Stein for pointing out this reference to them.

LEMMA 1. *Let $1 < p_1, p_2, q_1, q_2 \leq \infty$ and assume that $(p_1, p_2)$ and $(q_1, q_2)$ are good. Then*
$$\left( \frac{\theta}{p_1} + \frac{1-\theta}{q_1}, \frac{\theta}{p_2} + \frac{1-\theta}{q_2} \right)$$
*is good for all $0 < \theta < 1$.*

Next we need a duality lemma.

LEMMA 2. *Let $1 < p_1, p_2 < \infty$ such that $\frac{p_1 p_2}{p_1+p_2} \geq 1$. If $(p_1, p_2)$ is good, then so are the pairs*

$$\left(p_1, \left(\frac{p_1 p_2}{p_1 + p_2}\right)'\right) \quad \text{and} \quad \left(\left(\frac{p_1 p_2}{p_1 + p_2}\right)', p_2\right).$$

Here $p'$ denotes as usual the dual exponent of $p$. To prove the lemma, fix $\alpha \in \mathbb{R} \setminus \{0, -1\}$ and $f_1 \in \mathcal{S}(\mathbb{R})$ and consider the linear operator $H_\alpha(f_1, \_)$. The formal adjoint of this operator with respect to the natural bilinear pairing is

$$\text{sgn}(1 + \alpha) H_{-\frac{\alpha}{1+\alpha}}(f_1, \_),$$

as the following lines show:

$$\int \left(\text{p.v.} \int f_1(x-t) f_2(x+\alpha t) \frac{1}{t}\, dt\right) f_3(x)\, dx$$

$$= \text{p.v.} \int \int f_1(x - \alpha t - t) f_2(x) f_3(x - \alpha t)\, dx\, \frac{1}{t} dt$$

$$= \text{sgn}(1 + \alpha) \int \left(\text{p.v.} \int f_1(x-t) f_3(x - \frac{\alpha}{1+\alpha} t) \frac{1}{t}\, dt\right) f_2(x)\, dx.$$

Similarly, we observe that for fixed $f_2$ the formal adjoint of $H_\alpha(\_, f_2)$ is $-H_{-1-\alpha}(\_, f_2)$. This proves Lemma 2 by duality.

Now we are ready to prove estimate (2) in the remaining cases, i.e., for those pairs $(p_1, p_2)$ for which one of $p_1, p_2$ is smaller or equal two, and the other one is greater or equal two. In this case the constraint on $p$ is automatically satisfied. By symmetry it suffices to do this for $p_1 \in ]1, 2]$ and $p_2 \in [2, \infty]$. First observe that the pairs $(3,3)$ and $(3/2, 3/2)$ are good by the above propositions. Then the pairs $(2,2)$ and $(2, \infty)$ are good by interpolation and duality. Let $P$ be the set of all $p_1 \in ]1, 2]$ such that the pair $(p_1, p_2)$ is good for all $p_2 \in [2, \infty]$. The previous observations show that $2 \in P$. Define $p := \inf P$ and assume $p > 1$. Pick a small $\varepsilon > 0$ and a $p_1 \in P$ with $p_1 < p + \varepsilon$. If $\varepsilon$ is small enough, we can interpolate the good pairs $(p_1, \varepsilon^{-1})$ and $(1 + \varepsilon, 2 - \varepsilon)$ to obtain a good pair of the form $(q_\varepsilon, q_\varepsilon')$. Since $\lim_{\varepsilon \to 0} q_\varepsilon = \frac{3p-2}{2p-1} < p$ we have $q_\varepsilon < p$ provided $\varepsilon$ is small enough. By duality we see that the pair $(q, \infty)$ is good, and by Proposition 1 there is a $p_2 < 2$ such that $(q, p_2)$ is good. By interpolation $q \in P$ follows. This is a contradiction to $p = \inf P$; therefore the assumption $p > 1$ is false and we have $\inf P = 1$. Again by interpolation we observe $P = ]1, 2]$, which finishes the prove of estimate (2) for the remaining exponents.

### 3. Time-frequency decomposition of $H_\alpha$

In this section we write the bilinear operators $H_\alpha$ approximately as finite sums over rank one operators, each rank one operator being well localized in time and frequency. We mostly follow the corresponding section in [4], adopting the basic notation and definitions from there such as that of a phase plane representation.

In contrast to [4] we work out how the decomposition and the constants depend on $\alpha$, and we add an additional assumption (iv) in Proposition 3 which is necessary to prove $L^p$- estimates for $p < 2$. The reader should think of the functions $\theta_{\xi,\imath}$ in this assumption as being exponentials $\theta_{\xi,\imath}(x) = e^{i\eta_\imath x}$ for certain frequencies $\eta_\imath = \eta_\imath(\xi)$.

PROPOSITION 3. *Assume we are given exponents $1 < p_1, p_2 < 2$ such that $\frac{p_1 p_2}{p_1 + p_2} > \frac{2}{3}$, and we are given a constant $C_m$ for each integer $m \geq 0$. Then there is a constant $C$ depending on these data such that the following holds:*

*Let $S$ be a finite set, $\phi_1, \phi_2, \phi_3 : S \to \mathcal{S}(\mathbb{R})$ be injective maps, and $I, \omega_1, \omega_2, \omega_3 : S \mapsto \mathcal{J}$ be maps such that $I(S)$ is a grid, $\mathcal{J}_\omega := \omega_1(S) \cup \omega_2(S) \cup \omega_3(S)$ is a grid, and the following properties* (i)–(iv) *hold for all $\imath \in \{1, 2, 3\}$:*

(i) *The map*
$$\rho_\imath : \phi_\imath(S) \to \mathcal{R}, \ \phi_\imath(s) \mapsto I(s) \times \omega_\imath(s)$$
*is a phase plane representation with constants $C_m$.*

(ii) $\omega_\imath(s) \cap \omega_\jmath(s) = \emptyset$ *for all $s \in S$ and $\jmath \in \{1, 2, 3\}$ with $\imath \neq \jmath$.*

(iii) *If $\omega_\imath(s) \subset J$ and $\omega_\imath(s) \neq J$ for some $s \in S$, $J \in \mathcal{J}_\omega$, then $\omega_\jmath(s) \subset J$ for all $\jmath \in \{1, 2, 3\}$.*

(iv) *To each $\xi \in \mathbb{R}$ there is associated a measurable function $\theta_{\xi,\imath} : \mathbb{R} \to \{z \in \mathbb{C} : |z| = 1\}$ such that for all $s \in S$, $\jmath \in \{1, 2, 3\}$ and $J \in I(S)$ the following holds: If $\xi \in \omega_\jmath(s)$, $|J| \leq |I(s)|$, then*

(5) $$\inf_{\lambda \in \mathbb{C}} \|\phi_\imath(s) - \lambda \theta_{\xi,\imath}\|_{L^\infty(J)} \leq C_0 |J| \, |I(s)|^{-\frac{3}{2}} \left(1 + \frac{|c(J) - c(I(s))|}{|I(s)|}\right)^{-2}.$$

*For all $f_1, f_2 \in \mathcal{S}(\mathbb{R})$ and all maps $\varepsilon : S \to [-1, 1]$, we then have:*

(6) $$\left\| \sum_{s \in S} \varepsilon(s) |I(s)|^{-\frac{1}{2}} \langle f_1, \phi_1(s) \rangle \langle f_2, \phi_2(s) \rangle \phi_3(s) \right\|_{\frac{p_1 p_2}{p_1 + p_2}} \leq C \|f_1\|_{p_1} \|f_2\|_{p_2}.$$

In the rest of this section we prove Proposition 2 under the assumption that Proposition 3 above is true. Let $1 < p_1, p_2 < 2$ with $\frac{2}{3} < p := \frac{p_1 p_2}{p_1 + p_2}$ and $\alpha \in \mathbb{R} \setminus \{0, -1\}$.

Let $L$ be the smallest integer larger than
$$2^{10} \max\left\{|\alpha|, \frac{1}{|\alpha|}, \frac{1}{|1+\alpha|}\right\}.$$

The dependence on $\alpha$ will enter into our estimate via a polynomial dependence on $L$.

Define $\varepsilon := L^{-3}$. Pick a function $\psi \in \mathcal{S}(\mathbb{R})$ such that $\hat{\psi}$ is supported in $[L^3 - 1, L^3 + 1]$ and
$$\sum_{k\in\mathbb{Z}} \hat{\psi}(2^{\varepsilon k}\xi) = 1 \quad \text{for all} \quad \xi > 0.$$

Define
$$\psi_k(x) := 2^{-\frac{\varepsilon k}{2}} \psi(2^{-\varepsilon k}x)$$

and

(7) $$\tilde{H}_\alpha(f_1, f_2)(x) := \sum_{k\in\mathbb{Z}} 2^{-\frac{\varepsilon k}{2}} \int_{\mathbb{R}} f_1(x-t)f_2(x+\alpha t)\psi_k(t)dt.$$

It suffices to prove boundedness of $\tilde{H}_\alpha$. Pick a $\varphi \in \mathcal{S}(\mathbb{R})$ such that $\hat{\varphi}$ is supported in $[-1, 1]$ and

(8) $$\sum_{n,l\in\mathbb{Z}} \left\langle f, \varphi_{k,n,\frac{l}{2}} \right\rangle \varphi_{k,n,\frac{l}{2}} = f$$

for all Schwartz functions $f$, where
$$\varphi_{\kappa,n,l}(x) := 2^{-\frac{\varepsilon\kappa}{2}} \varphi(2^{-\varepsilon\kappa}x - n) e^{2\pi i 2^{-\varepsilon\kappa}xl}.$$

We apply this formula three times in (7) to obtain:

(9)
$$\tilde{H}_\alpha(f_1, f_2)(x) = \sum_{k,n_1,n_2,n_3,l_1,l_2,l_3\in\mathbb{Z}} C_{k,n_1,n_2,n_3,l_1,l_2,l_3} H_{k,n_1,n_2,n_3,l_1,l_2,l_3}(f_1, f_2)(x)$$

with
$$H_{k,n_1,n_2,n_3,l_1,l_2,l_3}(f_1, f_2)(x) := 2^{-\frac{\varepsilon k}{2}} \left\langle f_1, \varphi_{k,n_1,\frac{l_1}{2}} \right\rangle \left\langle f_2, \varphi_{k,n_2,\frac{l_2}{2}} \right\rangle \varphi_{k,n_3,\frac{l_3}{2}}(x)$$

and
(10)
$$C_{k,n_1,n_2,n_3,l_1,l_2,l_3} := \int\int \varphi_{k,n_1,\frac{l_1}{2}}(x-t)\varphi_{k,n_2,\frac{l_2}{2}}(x+\alpha t)\varphi_{k,n_3,\frac{l_3}{2}}(x)\psi_k(t)\,dt\,dx\ .$$

The proof of the following lemma is a straightforward calculation as in [4].

LEMMA 3. *There is a constant $C$ depending on $\phi$ and $\psi$ such that*

(11) $$|C_{k,n_1,n_2,n_3,l_1,l_2,l_3}| \leq C\left(1 + \frac{1}{L}\mathrm{diam}\{n_1, n_2, n_3\}\right)^{-100}.$$

*Moreover,*

$$C_{k,n_1,n_2,n_3,l_1,l_2,l_3} = 0,$$

*unless*

(12) $$l_1 \in \left[\left(-\frac{\alpha}{1+\alpha}l_3 + \frac{2}{1+\alpha}L^3\right) - L, \left(-\frac{\alpha}{1+\alpha}l_3 + \frac{2}{1+\alpha}L^3\right) + L\right]$$

*and*

(13) $$l_2 \in \left[\left(-\frac{1}{1+\alpha}l_3 - \frac{2}{1+\alpha}L^3\right) - L, \left(-\frac{1}{1+\alpha}l_3 - \frac{2}{1+\alpha}L^3\right) + L\right].$$

Now we can reduce Proposition 2 to the following lemma:

LEMMA 4. *There is a constant $C$ depending on $p_1$, $p_2$, $\varphi$, and $\psi$ such that the following holds:*

*Let $\nu > 0$ be an integer and let $S$ be a finite subset of $\mathbb{Z}^3$ such that for $(k,n,l), (k',n',l') \in S$ the following three properties are satisfied:*

(14) $\quad\quad\quad$ If $k \neq k'$, $\quad\quad$ then $|k - k'| > L^{10}$,

(15) $\quad\quad\quad$ if $n \neq n'$, $\quad\quad$ then $|n - n'| > L^{10}\nu$,

(16) $\quad\quad\quad$ if $l \neq l'$, $\quad\quad\quad$ then $|l - l'| > L^{10}$.

*Let $\nu_1$, $\nu_2$ be integers with $1 + \max\{|\nu_1|, |\nu_2|\} = \nu$ and let $\lambda_1, \lambda_2 : \mathbb{Z} \to \mathbb{Z}$ be functions such that $l_1 := \lambda_1(l_3)$ satisfies (12) and $l_2 := \lambda_2(l_3)$ satisfies (13) for all $l_3 \in \mathbb{Z}$. Then we have for all $f_1, f_2 \in \mathcal{S}(\mathbb{R})$ and all maps $\varepsilon : S \to [-1, 1]$:*

(17)
$$\left\| \sum_{(k,n,l) \in S} \varepsilon(k,n,l) H_{k,n+\nu_1,n+\nu_2,n,\lambda_1(l),\lambda_2(l),l}(f_1, f_2) \right\|_p \leq CL^{30}\nu^{10}\|f_1\|_{p_1}\|f_2\|_{p_2}.$$

Before proving the lemma we show how it implies boundedness of $\tilde{H}_\alpha$ and therefore proves Proposition 2. First observe that the lemma also holds without the finiteness condition on $S$. We can also remove conditions (14), (15), and (16) on $S$ at the cost of some additional powers of $L$ and $\nu$, so that the conclusion of the lemma without these hypotheses is

(18)
$$\left\| \sum_{(k,n,l) \in S} \varepsilon(k,n,l) H_{k,n+\nu_1,n+\nu_2,n,\lambda_1(l),\lambda_2(l),l}(f_1, f_2) \right\|_p \leq CL^{100}\nu^{20}\|f_1\|_{p_1}\|f_2\|_{p_2}.$$

Here we have used the quasi triangle inequality for $L^p$ which is uniform for $p > \frac{2}{3}$.

Observe that (18) and (11) imply

$$\text{(19)} \quad \left\| \sum_{(k,n,l) \in S} C_{k,n+\nu_1,n+\nu_2,n,\lambda_1(l),\lambda_2(l),l} H_{k,n+\nu_1,n+\nu_2,n,\lambda_1(l),\lambda_2(l),l}(f_1, f_2) \right\|_p$$
$$\leq C L^{200} \nu^{-50} \|f_1\|_{p_1} \|f_2\|_{p_2}.$$

Conditions (12) and (13) give a bound on the number of values the functions $\lambda_1$ and $\lambda_2$ can take at a fixed $l_3$ so that the coefficient $C_{k,n+\nu_1,n+\nu_2,n,\lambda_1(l),\lambda_2(l),l}$ does not vanish. Moreover there are of the order $\nu$ pairs $\nu_1, \nu_2$ such that $1 + \max\{|\nu_1|, |\nu_2|\} = \nu$. Hence,

$$\left\| \sum_{(k,n,l) \in S, n_1, n_2, l_1, l_2 \in \mathbb{Z}, 1+\max\{|n-n_1|,|n-n_2|\}=\nu} C_{k,n_1,n_2,n,l_1,l_2,l} H_{k,n_1,n_2,n,l_1,l_2,l}(f_1, f_2) \right\|_p$$
$$\leq C L^{300} \nu^{-20} \|f_1\|_{p_1} \|f_2\|_{p_2}.$$

Summing over all $\nu$ gives boundedness of $\tilde{H}_\alpha$.

It remains to prove Lemma 4. Clearly we intend to do this by applying Proposition 3. Fix data $S, \nu, \nu_1, \nu_2, \lambda_1, \lambda_2$ as in Lemma 4. Define functions $\phi_\imath : S \mapsto \mathcal{S}(\mathbb{R})$ as follows:

$$\phi_1(k, n, l) := L^{-10} \nu^{-2} \varphi_{k, n+\nu_1, \frac{\lambda_1(l)}{2}},$$
$$\phi_2(k, n, l) := L^{-10} \nu^{-2} \varphi_{k, n+\nu_2, \frac{\lambda_2(l)}{2}},$$
$$\phi_3(k, n, l) := L^{-10} \nu^{-2} \varphi_{k, n, \frac{l}{2}}.$$

If $E$ is a subset of $\mathbb{R}$ and $x \neq 0$ a real number we use the notation $x \cdot E := \{xy \in \mathbb{R} : y \in E\}$. This is not to be confused with the previously defined $xI$ for positive $x$ and intervals $I$. Pick three maps $\omega_1, \omega_2, \omega_3 : S \to \mathcal{J}$ such that the following properties (20)–(25) are satisfied for all $s = (k, n, l) \in S$:

$$\text{(20)} \quad -\frac{1+\alpha}{\alpha} \cdot \text{supp}\,(\widehat{\phi_1(s)}) \subset \omega_1(s),$$

$$\text{(21)} \quad -(1+\alpha) \cdot \text{supp}\,(\widehat{\phi_2(s)}) \subset \omega_2(s),$$

$$\text{(22)} \quad \text{supp}\,(\widehat{\phi_3(s)}) \subset \omega_3(s),$$

$$\text{(23)} \quad 2^{-\varepsilon(k+1)} L \leq |\omega_\imath(s)| \leq 2^{-\varepsilon k} L \text{ for } \imath = 1, 2, 3,$$

$$\text{(24)} \quad \mathcal{J}_\omega := \omega_1(S) \cup \omega_2(S) \cup \omega_3(S) \text{ is a grid},$$

and, for all $\imath, \jmath \in \{1, 2, 3\}$,

(25)        If $\omega_\imath(s) \subset J$ and $\omega_\imath(s) \neq J$ for some $J \in \mathcal{J}_\omega$, then $\omega_\jmath(s) \subset J$.

The existence of such a triple of maps is proved as in [4].

Next pick a map $I : S \to \mathcal{J}$ which satisfies the following three properties (26)–(28) for all $s = (k, n, l) \in S$:

(26)                $|c(I(s)) - 2^{\varepsilon k} n| \leq 2^{\varepsilon k} \nu,$

(27)                $2^4 2^{\varepsilon k} \nu \leq |I(s)| \leq 2^\varepsilon 2^4 2^{\varepsilon k} \nu,$

(28)                $I(S)$ is a grid.

The existence of such a map is again proved as in [4].

Now Lemma 4 follows immediately from the fact that the data $S$, $\phi_1$, $\phi_2$, $\phi_3$, $I$, $\omega_1$, $\omega_2$, and $\omega_3$ satisfy the hypotheses of Proposition 3. The verification of these hypotheses is as in [4] except for hypothesis (iv).

We prove hypothesis (iv) for $\imath = 1$, the other cases being similar. Define for $\xi \in \mathbb{R}$:
$$\theta_{\xi,1}(x) := e^{-2\pi i \frac{\alpha}{\alpha+1} x \xi}.$$

Pick $s = (k, n, l) \in S$. Obviously,
$$\nu^{-2} \varphi_{k, n+\nu_1, 0}(x) \leq C |I(s)|^{-\frac{1}{2}} \left(1 + \frac{|x - c(I(s))|}{|I(s)|}\right)^{-2}$$

and
$$\nu^{-2} (\varphi_{k, n+\nu_1, 0})'(x) \leq C |I(s)|^{-\frac{3}{2}} \left(1 + \frac{|x - c(I(s))|}{|I(s)|}\right)^{-2}.$$

Now let $\xi \in \omega_\jmath(s)$. By choice of $\theta_{\xi,1}$ we see that the function
$$\varphi_{k, n+\nu_1, \frac{\lambda_1(l)}{2}} \theta_{\xi,1}^{-1}$$
arises from $\varphi_{k, n+\nu_1, 0}$ by modulating with a frequency which is contained in $L^{10}[-|I(s)|^{-1}, |I(s)|^{-1}]$. Therefore,
$$(\phi_1(s) \theta_{\xi,1}^{-1})'(x) \leq C |I(s)|^{-\frac{3}{2}} \left(1 + \frac{|x - c(I(s))|}{|I(s)|}\right)^{-2}.$$

Now let $J \in I(S)$ with $|J| \leq |I(s)|$. Then we have
$$\inf_\lambda \|\phi_1(s) \theta_{\xi,1}^{-1} - \lambda\|_{L^\infty(J)} \leq |J| \left\|\left(\phi_1(s) \theta_{\xi,1}^{-1}\right)'\right\|_{L^\infty(J)}$$
$$\leq C |J| |I(s)|^{-\frac{3}{2}} \left(1 + \frac{|c(J) - c(I(s))|}{|I(s)|}\right)^{-2}.$$

This proves hypothesis (iv), and therefore finishes the reduction of Proposition 2 to Proposition 3.

## 4. Reduction to a symmetric statement

The following proposition is a variant of Proposition 3 which is symmetric in the indices 1, 2, and 3.

PROPOSITION 4. *Let $1 < p_1, p_2, p_3 < 2$ be exponents with*

$$1 < \frac{1}{p_1} + \frac{1}{p_2} + \frac{1}{p_3} < 2$$

*and let $C_m > 0$ for $m \geq 0$. Then there are constants $C, \lambda_0 > 0$ such that the following holds: Let $S, \phi_1, \phi_2, \phi_3, I, \omega_1, \omega_2, \omega_3$ be as in Proposition 3, let $f_i$, $i = 1, 2, 3$ be Schwartz functions with $\|f_i\|_{p_i} = 1$, and define*

$$E := \left\{ x \in \mathbb{R} : \max_i \left( M_{p_i}(Mf_i)(x) \right) \geq \lambda_0 \right\}.$$

*Then we have*

$$\sum_{s \in S : I(s) \not\subset E} |I(s)|^{-\frac{1}{2}} |\langle f_1, \phi_1(s) \rangle \langle f_2, \phi_2(s) \rangle \langle f_3, \phi_3(s) \rangle| \leq C.$$

We now prove that Proposition 3 follows from Proposition 4.
Let $1 < p_1, p_2 < 2$ and assume

$$p := \frac{p_1 p_2}{p_1 + p_2} > \frac{2}{3}.$$

Let $S, \phi_1, \phi_2, \phi_3, I, \omega_1, \omega_2, \omega_3, \varepsilon$ be as in the proposition and define for each $S' \subset S$

$$H_{S'}(f_1, f_2) = \sum_{s \in S'} \varepsilon(s) |I|^{-\frac{1}{2}} \langle f_1, \phi_1(s) \rangle \langle f_2, \phi_2(s) \rangle \phi_3(s).$$

By Marcinkiewicz interpolation ([2]), it suffices to prove a corresponding weak-type estimate instead of (6). By linearity and scaling invariance it suffices to prove that there is a constant $C$ such that for $\|f_1\|_{p_1} = \|f_2\|_{p_2} = 1$ we have

$$|\{x \in \mathbb{R} : |H_S(f_1, f_2)(x)| \geq 2\}| \leq C.$$

Pick an exponent $p_3$ such that the triple $p_1, p_2, p_3$ satisfies the conditions of Proposition 4, and let $\lambda_0$ be as in this proposition. Let $f_1$ and $f_2$ be Schwartz functions with $\|f_1\|_{p_1} = \|f_2\|_{p_2} = 1$.

Define

$$E_0 := \{x : \max\{M_{p_1}(Mf_1)(x), M_{p_2}(Mf_2)(x)\} \geq \lambda_0\}.$$

and

$$E_{\text{in}} := \left\{ x \in R : \left| H_{\{s \in S : I(s) \subset E_0\}}(f_1, f_2)(x) \right| \geq 1 \right\},$$
$$E_{\text{out}} := \left\{ x \in R : \left| H_{\{s \in S : I(s) \not\subset E_0\}}(f_1, f_2)(x) \right| \geq 1 \right\}.$$

It suffices to bound the measures of $E_{\text{in}}$ and $E_{\text{out}}$ by constants. We first estimate that of $E_{\text{out}}$ using Proposition 4. Let $\delta > 0$ be a small number and let $\theta : [0, \infty) \to [0, 1]$ be a smooth function which vanishes on the interval $[0, 1-\delta]$ and is constant equal to 1 on $[1, \infty)$. Extend this function to the complex plane by defining in polar coordinates $\theta(re^{i\phi}) := \theta(r)e^{-i\phi}$. Assume that $\delta$ is chosen sufficiently small to give

$$|E_{\text{out}}|^{\frac{1}{p_3}} \leq \left\| \theta\left(H_{\{s \in S: I(s) \not\subset E_0\}}(f_1, f_2)\right) \right\|_{p_3} \leq 2|E_{\text{out}}|^{\frac{1}{p_3}}.$$

Define

$$f_3 := \frac{\theta\left(H_{\{s \in S: I(s) \not\subset E_0\}}(f_1, f_2)\right)}{\left\| \theta\left(H_{\{s \in S: I(s) \not\subset E_0\}}(f_1, f_2)\right) \right\|_{p_3}}.$$

We can assume that $|E_{\text{out}}| > \lambda_0^{-p_3}$, because otherwise nothing is to prove. This assumption implies $\|M_{p_3}(Mf_3)\|_\infty < \lambda_0$. By applying Proposition 4, we obtain:

$$|E_{\text{out}}|^{1-\frac{1}{p_3}} \leq 2\left| \int H_{\{s \in S: I(s) \not\subset E_0\}}(f_1, f_2)(x) f_3(x)\, dx \right| \leq C.$$

Therefore $|E_{\text{out}}|$ is bounded by a constant.

It remains to estimate the measure of the set $E_{\text{in}}$, which is an elementary calculation. We need the following lemma:

LEMMA 5. *Let $J$ be an interval and define*

$$S_J := \{s \in S : I(s) = J\}.$$

*Then for all $m > 0$ there is a $C_m$ such that for all $A > 1$ and $f_1, f_2 \in \mathcal{S}(\mathbb{R})$ we have*:

$$\|H_{S_J}(f_1, f_2)\|_{L^1((AJ)^c)} \leq C_m |J| A^{-m} \left( \inf_{x \in J} M_{p_1} f_1(x) \right) \left( \inf_{x \in J} M_{p_2} f_2(x) \right).$$

We prove the lemma for $|J| = 1$, which suffices by homogeneity. For $m \geq 0$ define the weight

$$w_m(x) := (1 + \text{dist}(x, J))^m.$$

Then for $1 \leq r < 2$ we obtain the estimates

$$(29) \qquad \left\| \sum_{s \in S_J} \alpha_s \phi_\iota(s) \right\|_{L^{r'}(\omega_m)} \leq C_m \left\| (\alpha_s)_{s \in S_J} \right\|_{l^r(S_J)}$$

and

$$(30) \qquad \left\| (\langle f, \phi_\iota(s) \rangle)_{s \in S_J} \right\|_{l^{r'}(S_J)} \leq C_m \|f\|_{L^r(\omega_m^{-1})},$$

which follow easily by interpolation ([6]) from the trivial weighted estimate at $r = 1$ and the nonweighted estimate at $r = 2$.

Now define $r$ by
$$\frac{1}{r} = \frac{1}{p_1'} + \frac{1}{p_2'};$$
in particular we have $1 < r < 2$. By writing $H_{S_J}(f_1, f_2) = (H_{S_J}(f_1, f_2)w_m^{\frac{1}{r'}})w_m^{-\frac{1}{r'}}$ and applying Hölder we have for large $m$:
$$\|H_{S_J}(f_1, f_2)\|_{L^1((AJ)^c)} \leq C_M A^{-M} \|H_{S_J}(f_1, f_2)\|_{L^{r'}(w_m)}.$$

Here $M$ depends on $m$ and $r$ and can be made arbitrarily large by picking $m$ accordingly. By estimates (29) and (30) we can estimate the previously displayed expression further by
$$\leq C_M A^{-M} \left\| (\langle f_1, \phi_1(s) \rangle \langle f_2, \phi_2(s) \rangle)_{s \in S_J} \right\|_{l^r(S_J)}$$
$$\leq C_M A^{-M} \left\| (\langle f_1, \phi_1(s) \rangle)_{s \in S_J} \right\|_{l^{p_1'}(S_J)} \left\| (\langle f_2, \phi_2(s) \rangle)_{s \in S_J} \right\|_{l^{p_2'}(S_J)}$$
$$\leq C_M A^{-M} \|f_1\|_{L^{p_1}(w_{10}^{-1})} \|f_2\|_{L^{p_2}(w_{10}^{-1})}$$
$$\leq C_M A^{-M} \left( \inf_{x \in J} M_{p_1} f_1(x) \right) \left( \inf_{x \in J} M_{p_2} f_2(x) \right).$$

This finishes the proof of Lemma 5.

We return to the estimate of the set $E_{\text{in}}$. Define
$$E' := E_0 \cup \bigcup_{J \in I(S) : J \subset E} 4J.$$

Since $|E'| \leq 5|E_0| \leq C$, it suffices to prove

(31) $$\|H_{\{s \in S : I(s) \subset E_0\}}(f_1, f_2)\|_{L^1(E'^c)} \leq C.$$

Fix $k > 1$ and define
$$\mathcal{I}_k := \{J \in I(S) : J \subset E_0, 2^k J \subset E', 2^{k+1} J \not\subset E'\}.$$

Let $J \in \mathcal{I}_k$. Then for $\imath = 1, 2$ we have:
$$\inf_{x \in J} M_{p_\imath} f_\imath(x) \leq 2^{k+1} \inf_{x \in 2^{k+1}J} M_{p_\imath} f_\imath(x) \leq 2^{k+1},$$
since outside the set $E'$ the maximal function is bounded by 1. Hence, by the previous lemma,
$$\|H_{S_J}(f_1, f_2)\|_{L^1((E')^c)} \leq C_m |J| 2^{-km}.$$

Since $I(S)$ is a grid, it is easy to see that the intervals in $\mathcal{I}_k$ are pairwise disjoint; hence we have
$$\left\| H_{\{s \in S, I(s) \in \mathcal{I}_k\}}(f_1, f_2) \right\|_{L^1((E')^c)} \leq C_m |E_0| 2^{-km}.$$

By summing over all $k > 1$ we prove (31). This finishes the estimate of the set $|E_{\text{in}}|$ and therefore the reduction of Proposition 3 to Proposition 4.

## 5. The combinatorics on the set $S$

We prove Proposition 4. Let $1 < p_1, p_2, p_3 < 2$ be exponents with
$$1 < \frac{1}{p_1} + \frac{1}{p_2} + \frac{1}{p_3} < 2.$$
Let $\eta > 0$ be the largest number such that $\frac{1}{\eta}$ is an integer and
$$\eta \leq 2^{-100} \left(2 - \sum_i \frac{1}{p_i}\right) \min_j \left(1 - \frac{1}{p_j}\right).$$

Let $S$, $\phi_1$, $\phi_2$, $\phi_3$, $I$, $\omega_1$, $\omega_2$, and $\omega_3$ be as in Propositions 3 and 4. Let $f_i$, $i = 1, 2, 3$ be Schwartz functions with $\|f_i\|_{p_i} = 1$. Without loss of generality we can assume that for all $s \in S$,

(32) $$I(s) \not\subset \left\{x \in \mathbb{R} : \max_i (M_{p_i}(Mf_i)(x)) \geq \lambda_0\right\},$$

where $\lambda_0$ is a constant which we will specify later.

Define a partial order $\ll$ on the set of rectangles by

(33) $$J_1 \times J_2 \ll J_1' \times J_2' \ , \ \text{if } J_1 \subset J_1' \text{ and } J_2' \subset J_2.$$

A subset $T \subset S$ is called a *tree of type $i$*, if the set $\rho_i(T)$ has exactly one maximal element with respect to $\ll$. This maximal element is called the *base* of the tree $T$ and is denoted by $s_T$. Define $J_T := I(s_T)$.

Define $S_{-1} := S$. Let $k \geq 0$ be an integer and assume by recursion that we have already defined $S_{k-1}$. Define
$$S_k := S_{k-1} \setminus \bigcup_{i,j=1}^{3} \left(\bigcup_{l=0}^{\infty} T_{k,i,j,l}\right),$$
where the sets $T_{k,i,j,l}$ are defined as follows. Let $k \geq 0$ and $i, j \in \{1, 2, 3\}$ be fixed. Let $l \geq 0$ be an integer and assume by recursion that we have already defined $T_{k,i,j,\lambda}$ for all integers $\lambda$ with $0 \leq \lambda < l$. If one of the sets $T_{k,i,j,\lambda}$ with $\lambda < l$ is empty, then define $T_{k,i,j,l} := \emptyset$. Otherwise let $\mathcal{F}$ denote the set of all trees $T$ of type $i$ which satisfy the following conditions (34)–(36):

(34) $$T \subset S_{k-1} \setminus \bigcup_{\lambda < l} T_{k,i,j,\lambda},$$

(35) if $i = j$, then $|\langle f_j, \phi_j(s) \rangle| \geq 2^{-\eta k} 2^{-\frac{k}{p_j'}} |I(s)|^{\frac{1}{2}}$ for all $s \in T$,

(36) if $i \neq j$, then $\left\|\left(\sum_{s \in T} \frac{|\langle f_j, \phi_j(s) \rangle|^2}{|I(s)|} 1_{I(s)}\right)^{\frac{1}{2}}\right\|_1 \geq 2^4 2^{-\frac{k}{p_j'}} |J_T|.$

If $\mathcal{F}$ is empty, then we define $T_{k,i,j,l} := \emptyset$. Otherwise define $\mathcal{F}_{\max}$ to be the set of all $T_{\max} \in \mathcal{F}$ which satisfy:

(37) if $T \in \mathcal{F}$ , $T_{\max} \subset T$ , then $T = T_{\max}$.

Choose $T_{k,\imath,\jmath,l} \in \mathcal{F}_{\max}$ such that for all $T \in \mathcal{F}_{\max}$,

(38) $\quad\quad\quad$ if $\imath < \jmath$, $\quad\quad$ then $\omega_\imath(s_{T_{k,\imath,\jmath,l}}) \not< \omega_\imath(s_T)$,

(39) $\quad\quad\quad$ if $\imath > \jmath$, $\quad\quad$ then $\omega_\imath(s_T) \not< \omega(s_{T_{k,\imath,\jmath,l}})$.

Here $[a,b[ \not< [a',b'[$ means $b > a'$. Observe that $T_{k,\imath,\jmath,l}$ actually satisfies (38) and (39) for all $T \in \mathcal{F}$. This finishes the definition of the sets $T_{k,\imath,\jmath,l}$ and $S_k$.

Since $S$ is finite, $T_{k,\imath,\jmath,l} = \emptyset$ for sufficiently large $l$. In particular, each $s \in S_k$ satisfies

(40) $$|\langle f_\imath, \phi_\imath(s)\rangle| \leq 2^{-\eta k} 2^{-\frac{k}{p_\imath'}} |I(s)|^{\frac{1}{2}}$$

for all $\imath$, since the set $\{s\}$ is a tree of type $\imath$ which by construction of $S_k$ does not satisfy (35) for $\jmath = \imath$. Similarly for $j \neq i$ each tree $T \subset S_k$ of type $\imath$ satisfies

(41) $$\left\| \left( \sum_{s \in T} \frac{|\langle f_\jmath, \phi_\jmath(s)\rangle|^2}{|I(s)|} 1_{I(s)} \right)^{\frac{1}{2}} \right\|_1 \leq 2^4 2^{-\frac{k}{p_\jmath'}} |J_T|.$$

Moreover, (40) implies that the intersection of all $S_k$ contains only elements $s$ with $\prod_j \langle f_\jmath, \phi_\jmath(s)\rangle = 0$.

Let $k \leq \eta^{-2}$ and assume $T_{k,\imath,\jmath,l}$ is a tree. Observe that (35) and (36) together with Lemma 6 in Section 7 provide a lower bound on the maximal function $M_{p_j}(Mf_j)(x)$ for $x \in J_{T_{k,\imath,\jmath,l}}$. This lower bound depends only on $\eta$, $p_\jmath$ and the constants $C_m$ of the phase plane representation. Therefore if we choose the constant $\lambda_0$ in (32) small enough depending on $\eta$, $p_\jmath$, and $C_m$, it then is clear that $T_{k,\imath,\jmath,l} = \emptyset$ for $k \leq \eta^{-2}$.

Now we have

$$\sum_{s \in S} |I(s)|^{-\frac{1}{2}} \prod_j |\langle f_\jmath, \phi_\jmath(s)\rangle| \leq \sum_{k > \eta^{-2}} \sum_{\imath,\jmath} \sum_{l=0}^{\infty} \left( \sup_{s \in T_{k,\imath,\jmath,l}} |I(s)|^{-\frac{1}{2}} |\langle f_\imath, \phi_\imath(s)\rangle| \right)$$
$$\times \prod_{\kappa \neq \imath} \left( \sum_{s \in T_{k,\imath,\jmath,l}} |\langle f_\kappa, \phi_\kappa(s)\rangle|^2 \right)^{\frac{1}{2}}.$$

Using (40), (41) and Lemma 7 of Section 7 we can bound this by

$$\leq C \sum_{k > \eta^{-2}} 2^{-\sum_\jmath \frac{k}{p_\jmath'}} \sum_{\imath,\jmath} \sum_{l=0}^{\infty} |J_{T_{k,\imath,\jmath,l}}|.$$

Now we apply the estimate

(42) $$\sum_{l=0}^{\infty} |J_{T_{k,\imath,\jmath,l}}| \leq C 2^{10\eta p_\jmath' k} 2^k,$$

for each $k > \eta^{-2}, \imath\jmath$, which is proved in Sections 6 and 8. This bounds the previously displayed expression by

$$\text{(43)} \qquad \leq C \sum_{k > \eta^{-2}} 2^{-\sum \frac{k}{p_j'}} 2^{10\eta p_j' k} 2^k.$$

This is less than a constant since

$$\sum_j \frac{1}{p_j'} \geq 1 + 10\eta \max_j p_j'$$

by the choice of $\eta$. This finishes the proof of Proposition 4 up to the proof of estimate (42) and Lemmata 6 and 7.

## 6. Counting the trees for $\imath = \jmath$

We prove estimate (42) in the case $\imath = \jmath$. Thus fix $k > \eta^{-2}, \imath, \jmath$ with $\imath = \jmath$. Let $\mathcal{F}$ denote the set of all trees $T_{k,\imath,\jmath,l}$. Observe that for $T, T' \in \mathcal{F}, T \neq T'$ we have, by (37), that $T \cup T'$ is not a tree; therefore

$$\rho_\imath(s_T) \cap \rho_\imath(s_{T'}) = \emptyset.$$

Define $b := 2^{-\eta k} 2^{-\frac{k}{p_\imath'}}$. Then by (35) for all $T \in \mathcal{F}$

$$\text{(44)} \qquad |\langle f_\imath, \phi_\imath(s_T)\rangle| \geq b |J_T|^{\frac{1}{2}}.$$

Finally recall that for all $s \in S$:

$$\text{(45)} \qquad I(s) \not\subset \{x : M_{p_\imath}(Mf_\imath)(x) \geq \lambda_0\}.$$

Our proof goes in the following four steps:

*Step* 1. Define the counting function

$$\text{(46)} \qquad N_\mathcal{F}(x) := \sum_{T \in \mathcal{F}} 1_{J_T}(x).$$

We have to estimate the $L^1$-norm of the counting function. Since the counting function is integer-valued, it suffices to show a weak-type $1 + \varepsilon$ estimate for small $\varepsilon$. More precisely it suffices to show for all integers $\lambda \geq 1$ and sufficiently small $\delta, \varepsilon > 0$, $\delta = \delta(\eta, p_\imath), \varepsilon = \varepsilon(\eta, p_\imath)$:

$$|\{x \in \mathbb{R} : N_\mathcal{F}(x) \geq \lambda\}| \leq b^{-p_\imath' - \delta} \lambda^{-1-\varepsilon}.$$

Fix such a $\lambda$. As in [4] there is a subset $\mathcal{F}' \subset \mathcal{F}$ such that, if we define $N_{\mathcal{F}'}$ analogously to $N_\mathcal{F}$,

$$\{x \in \mathbb{R} : N_{\mathcal{F}'}(x) \geq \lambda\} = \{x \in \mathbb{R} : N_\mathcal{F}(x) \geq \lambda\}$$

and $\|N_{\mathcal{F}'}\|_\infty \leq \lambda$. This is due to the grid structure of $I(S)$.

*Step* 2. Let $A > 1$ be a number whose value will be specified later. We can write

$$\mathcal{F}' = \left(\bigcup_{m=1}^{A^{10}} \mathcal{F}_m\right) \cup \mathcal{F}'' \tag{47}$$

such that if $T, T' \in \mathcal{F}_m$ for some $m$ and $T \neq T'$, then

$$(AJ_T \times \omega(s_T)) \cap (AJ_{T'} \times \omega(s_{T'})) = \emptyset,$$

and

$$\sum_{T \in \mathcal{F}''} |J_T| \leq Ce^{-A} \sum_{T \in \mathcal{F}_1} |J_T|. \tag{48}$$

For a proof of this fact see the proof of the separation lemma in [4].

*Step* 3. Let $1 \leq m \leq A^{10}$. The following lines hold for all sufficiently small $\delta, \varepsilon > 0$. The arguments may require $\delta, \varepsilon$ to change from line to line. For a tempered distribution $f$, $x \in \mathbb{R}$, and $T \in \mathcal{F}_m$ define

$$Bf(x)(T) := \frac{\langle f, \phi_\iota(s_T)\rangle}{|J_T|^{\frac{1}{2}}} 1_{J_T}(x).$$

Let $L^2(\mathbb{R}, l^2(\mathcal{F}))$ be the Banach space of square-integrable functions on $\mathbb{R}$ with values in $l^2(\mathcal{F})$, and analogously for other exponents. Then we have the following estimate by Lemma 4.3 in [4]

$$\|Bf\|_{L^2(\mathbb{R}, l^2(\mathcal{F}_m))} = \left(\sum_{T \in \mathcal{F}_m} |\langle f_\iota, \phi_\iota(s_T)\rangle|^2\right)^{\frac{1}{2}} \leq C(1 + A^{-\frac{1}{\varepsilon}}\lambda)\|f\|_2.$$

We also trivially have

$$\|Bf\|_{L^{1+\delta}(\mathbb{R}, l^\infty(\mathcal{F}_m))} = \left(\int \left(\sup_{T \in \mathcal{F}_m : x \in J_T} \frac{|<f_\iota, \phi_\iota(s_T)>|}{|J_T|^{\frac{1}{2}}}\right)^{1+\delta} dx\right)^{\frac{1}{1+\delta}}$$

$$\leq C\|Mf\|_{1+\delta} \leq C\|f\|_{1+\delta}.$$

By interpolation we have for $1 < p < 2$:

$$\|Bf\|_{L^p(\mathbb{R}, l^{p'+\delta}(\mathcal{F}_m))} \leq C(1 + A^{-\frac{1}{\varepsilon}}\lambda)\|f\|_p.$$

Let $J \in I(S)$, and let $\mathcal{F}_{m,J}$ be the set of $T \in \mathcal{F}_m$ such that $J_T \subset J$. By a localization argument, as in [4], we see that

$$\|Bf\|_{L^p(\mathbb{R}, l^{p'+\delta}(\mathcal{F}_{m,J}))} \leq C\lambda^\varepsilon(1 + A^{-\frac{1}{\varepsilon}}\lambda)|J|^{\frac{1}{p}} \inf_{x \in J} M_p(Mf)(x).$$

In the following, $g^\sharp$ denotes the sharp maximal function of $g$ with respect to the given grid, as in [4]. We define $N_{\mathcal{F}_m,J}$ in analogy to (46) to be the counting function of the trees $T \in \mathcal{F}_m$ for which $I_T \subset J$. We apply the previous estimate for $f_\imath$ and use (44) to obtain

$$\left(N_{\mathcal{F}_m}^{\frac{p}{p'+\delta}}\right)^\sharp(x) \leq \sup_{J:x\in J}\left(\frac{1}{|J|}\int_J N_{\mathcal{F}_m,J}(x)^{\frac{p}{p'+\delta}}\,dx\right)$$

$$\leq b^{-p}\sup_{J:x\in J}\frac{1}{|J|}\left\|\left(\sum_{T\in\mathcal{F}_{m,J}}\frac{|\langle f_\imath,\phi_\imath(s_T)\rangle|^{p'+\delta}}{|J_T|^{\frac{p'+\delta}{2}}}1_{J_T}\right)^{\frac{1}{p'+\delta}}\right\|_p^p$$

$$\leq b^{-p}C\left(\lambda^\varepsilon(1+A^{-\frac{1}{\varepsilon}}\lambda)M_p(Mf_\imath)(x)\right)^p.$$

Using (45) we can sharpen this argument in the case $p = p_\imath$ to

$$\left(N_{\mathcal{F}_m}^{\frac{p_\imath}{p_\imath'+\delta}}\right)^\sharp(x) \leq Cb^{-p_\imath}\left(\lambda^\varepsilon(1+A^{-\frac{1}{\varepsilon}}\lambda)\min\{M_{p_\imath}(Mf_\imath)(x),\lambda_0\}\right)^{p_\imath}.$$

Taking the $\frac{p_\imath'+2\delta}{p_\imath}$- norm on both sides and raising to the $\frac{p_\imath'+\delta}{p_\imath}$-th power gives

(49) $$\|N_{\mathcal{F}_m}\|_{\frac{p_\imath'+2\delta}{p_\imath'+\delta}} \leq Cb^{-p_\imath'-\delta}\left(\lambda^\varepsilon(1+A^{-\frac{1}{\varepsilon}}\lambda)\right)^{p_\imath'+\delta}.$$

*Step* 4. We split the counting function $N_{\mathcal{F}'}$ according to (47) and use the weak-type estimate following from (49) on the first part and estimate (48) together with (49) and the fact that the counting function is integer-valued on the second part. This gives

$$\left\{x\in\mathbb{R}: N_{\mathcal{F}'}(x)\geq A^{10}\lambda\right\} \leq CA^{10}\lambda^{-\frac{p_\imath'+2\delta}{p_\imath'+\delta}}b^{-p_\imath'-2\delta}\left(\lambda^\varepsilon(1+A^{-\frac{1}{\varepsilon}}\lambda)\right)^{p_\imath'+2\delta}$$
$$+e^{-A}Cb^{-p_\imath'-2\delta}\left(\lambda^\varepsilon(1+A^{-\frac{1}{\varepsilon}}\lambda)\right)^{p_\imath'+2\delta}.$$

Choosing $A$ of the order $\lambda^\varepsilon$ and $\varepsilon \ll \delta$ gives

$$\{x\in\mathbb{R}: N_{\mathcal{F}'}(x)\geq \lambda\} \leq C\lambda^{-1-\varepsilon}b^{-p_\imath'-\delta}.$$

According to Step 1 this finishes the proof of estimate (42) in the case $\imath = \jmath$.

## 7. Estimates on a single tree

This section collects some standard facts from Calderón-Zygmund theory, adapted to the setup of trees.

LEMMA 6. *Fix* $k, \imath, \jmath, l$ *such that* $T := T_{k,\imath,\jmath,l}$ *is a tree, assume* $\imath \neq \jmath$, *and let* $1 < p \leq 2$. *We then have*

(50) $$\left\|\left(\sum_{s\in T}\frac{|\langle f,\phi_\jmath(s)\rangle|^2}{|I(s)|}1_{I(s)}\right)^{\frac{1}{2}}\right\|_p \leq C\|f\|_p.$$

*For each interval $J \in I(S)$ define $T_J := \{s \in T : I(s) \subset J\}$. Then we obtain*

$$\left\|\left(\sum_{s \in T_J} \frac{|\langle f, \phi_J(s) \rangle|^2}{|I(s)|} 1_{I(s)}\right)^{\frac{1}{2}}\right\|_p \leq C |J|^{\frac{1}{p}} \inf_{x \in J} M_p(Mf)(x). \tag{51}$$

*For each $s \in T$, let $h_s$ be a measurable function supported in $I(s)$ with $\|h_I(x)\|_\infty = |I(s)|^{-\frac{1}{2}}$, $\|h\|_2 = 1$, and $\langle h_s, h_{s'} \rangle = 0$ for $s \neq s'$. Then for all maps $\varepsilon : T \to \{-1, 1\}$, we have*

$$\left\|\sum_{s \in T} \varepsilon(s) \langle f, \phi_J(s) \rangle h_s\right\|_p \leq C \|f\|_p. \tag{52}$$

First we prove estimate (52). The estimate is true in the case $p = 2$, as is proved in [4]. By interpolation it suffices to prove the weak-type estimate

$$\left|\left\{x \in \mathbb{R} : \sum_{s \in T} \varepsilon(s) \langle f, \phi_J(s) \rangle h_s(x) \geq C\lambda\right\}\right| \leq C' \frac{\|f\|_1}{\lambda}. \tag{53}$$

Let $f \in L^1(\mathbb{R})$. We write $f$ as the sum of a good function $g$ and a bad function $b$ as follows. Let $\{I_n\}_n$ be the set of maximal intervals of the grid $I(S)$ for which

$$\int_{I_n} |f(x)| \, dx \geq \lambda |I_n|.$$

Let $\xi \in \omega_\imath(s_T)$, and pick a function $\theta_{\xi,\imath}$ as in hypothesis (iv) of Proposition 3. For each of the intervals $I_n$, define

$$b_n(x) := 1_{I_n}(x) \left( f(x) - \lambda_n \theta_{\xi,\imath}(x) \right),$$

where $\lambda_n$ is chosen such that $b_n$ is orthogonal to $\theta_{\xi,\imath}$. Obviously $\lambda_n$ is bounded by $C\|f(x)\|_{L^1(I_n)}$. Define $b := \sum_n b_n$ and $g := f - b$. It suffices to prove estimate (53) for the good and bad function separately. The estimate for the good function follows immediately from estimate (52) for $p = 2$. For the bad function we proceed as follows. Since the set

$$E := \bigcup_n 2 I_n$$

is bounded in measure by $C\lambda^{-1}$, it suffices to prove the strong-type estimate

$$\left\|\sum_n \left(\sum_{s \in T} \varepsilon(s) \langle b_n, \phi_J(s) \rangle h_s\right)\right\|_{L^1(E^c)} \leq C \|f\|_1. \tag{54}$$

We estimate each summand separately. Obviously, we have

$$\left\|\sum_{s \in T} \varepsilon(s) \langle b_n, \phi_J(s) \rangle h_s\right\|_{L^1(E^c)} \leq \sum_{s \in T : I(s) \not\subset 2 I_n} |I(s)|^{\frac{1}{2}} |\langle b_n, \phi_J(s) \rangle|.$$

For each integer $k$ let $T_k$ be the set of those $s \in T$, for which $|I(s)| \leq 2^k |I_n| < 2|I(s)|$ and $I(s) \not\subset 2I_n$. For $k < 2$ we use the estimate

$$\sum_{s \in T_k} |I(s)|^{\frac{1}{2}} |\langle b_n, \phi_J(s) \rangle| \tag{55}$$

$$\leq C \|b_n\|_1 \sum_{s \in T_k} \left(1 + \frac{|c(I(s)) - c(I_n)|}{|I(s)|}\right)^{-2}$$

$$\leq C \|b_n\|_1 \int_{(2I_n)^c} \sum_{s \in T_k} \frac{1}{2^k |I_n|} \left(1 + \frac{x - c(I_n)|}{2^k |I_n|}\right)^{-2} 1_{I(s)}(x)\, dx$$

$$\leq C \|b_n\|_1 2^k.$$

For the last inequality we have seen that the intervals $I(s)$ with $s \in T_k$ are pairwise disjoint.

For $k > 2$ we use the orthogonality of $b_n$ and $\theta_{\xi,\iota}$ as well as hypothesis (iv) of Proposition 3 to obtain

$$\sum_{s \in T_k} |I(s)|^{\frac{1}{2}} |\langle b_n, \phi_J(s) \rangle| \leq \sum_{s \in T_k} |I(s)|^{\frac{1}{2}} \|b_n\|_1 \inf_\lambda \|\phi_J(s) - \lambda \theta_{\xi,\iota}\|_{L^\infty(I_n)} \tag{56}$$

$$\leq C \|b_n\|_1 \sum_{s \in T_k} \left(1 + \frac{|c(I(s)) - c(I_n)|}{|I(s)|}\right)^{-2} \frac{|I_n|}{|I(s)|}$$

$$\leq C \|b_n\|_1 2^{-k}.$$

The last inequality follows by a similar argument as in the case $k \leq 2$. Summing (55) and (56) over $k$ and $n$ gives (54) and finishes the proof of (52).

We prove estimate (50). Observe that (52) is not void, since functions $h_s$ clearly exist. Therefore we can average (52) over all choices of $\varepsilon$ to obtain:

$$2^{-|T|} \sum_\varepsilon \left\| \sum_{s \in T} \varepsilon(s) \langle f, \phi_J(s) \rangle h_s \right\|_p^p = \int_\mathbb{R} 2^{-n} \sum_\varepsilon \left( \sum_{s \in T} \varepsilon(s) \langle f, \phi_J(s) \rangle h_s(x) \right)^p dx$$

$$\leq C \|f\|_p^p.$$

Now Khinchine's inequality gives

$$\int_\mathbb{R} \left( 2^{-n} \sum_\varepsilon \left( \sum_{s \in T} \varepsilon(s) \langle f, \phi_J(s) \rangle h_s(x) \right)^2 \right)^{\frac{p}{2}} dx \leq C \|f\|_p^p,$$

which immediately implies estimate (50).

To prove (51) fix a $J$ and write $f = f 1_{2J} + f 1_{(2J)^c}$. It suffices to prove the estimate separately for both summands. For the first summand we simply apply (50). For the second summand we write

$$\left(\sum_{s \in T_J} \frac{|\langle f 1_{(2J)^c}, \phi_J(s)\rangle|^2}{|I(s)|} 1_{I(s)}(x)\right)^{\frac{1}{2}} \leq C \sum_{s \in T_J : x \in I(s)} Mf(x)|I(s)||J|^{-1}$$
$$\leq CMf(x) 1_J(x).$$

The last inequality follows by summing a geometric series. This proves the estimate for the second summand and finishes the proof of Lemma 6.

LEMMA 7. *Fix $k \geq \eta^{-2}, \imath, \jmath, l$ such that $T := T_{k,\imath,\jmath,l}$ is a tree and assume $\imath \neq \jmath$. Then we have*

$$(57) \quad \left(\sum_{s \in T} |\langle f_\jmath, \phi_\jmath(s)\rangle|^2\right)^{\frac{1}{2}} \leq C \left\|\left(\sum_{s \in T} \frac{|\langle f_\jmath, \phi_\jmath(s)\rangle|^2}{|I(s)|} 1_{I(s)}\right)^{\frac{1}{2}}\right\|_1 |J_T|^{-\frac{1}{2}}.$$

*Proof.* For each $J \in I(S)$,

$$(58) \quad \frac{1}{|J|} \int_J \left(\left(\sum_{s \in T : I(s) \subset J} \frac{|\langle f_\jmath, \phi_\jmath(s)\rangle|^2}{|I(s)|} 1_{I(s)}(x)\right)^{\frac{1}{2}}\right) dx \leq C 2^{-\frac{k}{p_{\jmath'}}},$$

since the set $\{s \in T : I(s) \subset J\}$ is a union of trees $\{T_n\}_n$ which satisfy (41) for $k-1$ and

$$\sum_n |J_{T_n}| \leq |J_T|.$$

Define for $x \in \mathbb{R}$ and $s \in T$:

$$F(x)(s) := \sum_{s \in T} \frac{|\langle f_\jmath, \phi_\jmath(s)\rangle|^2}{|I(s)|} 1_{I(s)}(x).$$

Since $F$ is supported on $J_T$, we have

$$\|F\|_{L^2(\mathbb{R}, l^2(T))} \leq |J_T|^{\frac{1}{2}} \|F\|_{\text{BMO}(\mathbb{R}, l^2(T))}.$$

Here BMO is understood with respect to the grid $I(S)$ as in [4]. We prove Lemma 7 by estimating this BMO-norm with (58) and (36).

## 8. Counting the trees for $\imath \neq \jmath$

We prove estimate (42) in the case $\imath \neq \jmath$. Thus fix $k \geq \eta^{-2}, \imath, \jmath$ with $\imath \neq \jmath$. Let $\mathcal{F}$ denote the set of all trees $T_{k,\imath,\jmath,l}$.

As in [4] we define for $T \in \mathcal{F}$:

$$
\begin{aligned}
T^{\min} &:= \{s \in T : \rho_\iota(s) \text{ is minimal in } \rho_\iota(T)\}, \\
T^{\text{fat}} &:= \{s \in T : 2^5 2^{\eta k} |I(s)| \geq |J_T|\}, \\
T^\partial &:= \{s \in T : I(s) \cap (1 - 2^{-4}) J_T = \emptyset\}, \\
T^{\partial \max} &:= \{s \in T^\partial : \rho_\iota(s) \text{ is maximal in } \rho_\iota(T^\partial)\}, \\
T^{\text{nice}} &:= T \setminus \left(T^{\min} \cup T^{\text{fat}} \cup T^\partial\right).
\end{aligned}
$$

Define $b := 2^{-\frac{k}{p_J'}}$. By similar arguments as in [4] we have the estimate

(59) \quad if $\iota \neq \jmath$, then $\left\| \left( \sum_{s \in T^{\text{nice}}} \frac{|\langle f_\jmath, \phi_\jmath(s) \rangle|^2}{|I(s)|} 1_{I(s)} \right)^{\frac{1}{2}} \right\|_1 \geq b |J_T|.$

Define the counting function

$$N_{\mathcal{F}}(x) := \sum_{T \in \mathcal{F}} 1_{J_T}(x).$$

As in Section 6 it suffices to show

(60) \quad $|\{x \in \mathbb{R} : N_{\mathcal{F}'}(x) \geq \lambda\}| \leq b^{-p_J' - \delta} \lambda^{-1-\varepsilon}$

for all integers $\lambda \geq 1$ and small $\varepsilon, \delta > 0$. In addition, we can assume that $\|N_{\mathcal{F}}\|_\infty \leq \lambda$.

Let $y \in \mathbb{R}$, $T \in \mathcal{F}$, $x \in J_T$, and $s \in T$. For $f \in \mathcal{S}(\mathbb{R})$ define

$$Sf(y)(T)(x)(s) := \frac{\langle f, \phi_\jmath(s) \rangle}{|I(s)|^{\frac{1}{2}}} 1_{I(s)}(x) 1_{J_T}(y).$$

Consider $J_T$ as a measure space with Lebesgue measure normalized to 1. Then the operator is bounded from $L^2$ to $L^2(\mathbb{R}, l^2(\mathcal{F}, (L^2(J_T, l^2(T)))))$, as we see below. We have used a sloppy notation for the second Banach space: The range space $L^2(J_T, l^2(T))$ depends on the variable $T \in \mathcal{F}$. To make this space independent of $T$, we take the direct sum of these Banach spaces as $T$ varies over $\mathcal{F}$, and we let $Sf(y)(T)$ be nonzero only on the component corresponding to $T$. This is how we interpret the above notation. To see the claimed estimate we calculate:

$$
\int \sum_{T \in \mathcal{F}} \frac{1}{|J_T|} \int \sum_{s \in T} \frac{|\langle f, \phi_\jmath(s) \rangle|^2}{|I(s)|} 1_{I(s)}(x) 1_{J_T}(y) \, dx \, dy = \sum_{s \in \bigcup_{T \in \mathcal{F}} T} |\langle f, \phi_\jmath(s) \rangle|^2
$$
$$\leq C(1 + \lambda A^{-\frac{1}{\varepsilon}}) \|f\|_2^2,$$

the last inequality being taken from [4]. The operator is also bounded from $L^{1+2\delta}$ into

$$L^{1+2\delta}(\mathbb{R}, l^\infty(\mathcal{F}, L^{1+\delta}(J_T, l^2(T))))$$

since by Lemma 6:

$$\int \left( \sup_{T \in \mathcal{F}} \left( \frac{1}{|J_T|} \int \left( \sum_{s \in T} \left( \frac{|\langle f, \phi_J(s) \rangle|}{|I(s)|^{\frac{1}{2}}} 1_{I(s)}(x) 1_{J_T}(y) \right)^2 \right)^{\frac{1+\delta}{2}} dx \right)^{\frac{1}{1+\delta}} \right)^{1+2\delta} dy$$

$$\leq \int \left( \sup_{T \in \mathcal{F}: y \in J_T} \frac{1}{|J_T|} \left\| \left( \sum_{s \in T} \left( \frac{|\langle f, \phi_J(s) \rangle|}{|I(s)|^{\frac{1}{2}}} 1_{I(s)} \right)^2 \right)^{\frac{1}{2}} \right\|_{1+\delta}^{1+\delta} \right)^{\frac{1+2\delta}{1+\delta}} dy$$

$$\leq C \int (M_{1+\delta}(Mf)(y))^{1+2\delta} \, dy$$

$$\leq C \|f\|_{1+2\delta}^{1+2\delta}.$$

By complex interpolation and the fact that $L^q(J_T) \subset L^1(J_T)$ for $q \geq 1$ we obtain that $S$ maps $L^p$ into $L^p(\mathbb{R}, l^{p'+\delta}(\mathcal{F}, L^1(J_T, l^2(T))))$ with norm less than $C(1 + \lambda A^{-\frac{1}{\varepsilon}})$.

Let $J \in I(S)$ and define $\mathcal{F}_J$ to be the set of $T \in \mathcal{F}$ such that $J_T \subset J$. Then we can localize as before to get

$$\|Sf\|_{L^p(\mathbb{R}, l^{p'+\delta}(\mathcal{F}_J, L^1(J_T, l^2(T))))} \leq C\lambda^\varepsilon (1 + \lambda A^{-\frac{1}{\varepsilon}})|J|^{\frac{1}{p}} \inf_{x \in J} M^p(Mf)(x).$$

Using the estimate (59) on nice trees gives, for $f = f_J$ and $p = p_J$,

$$\left( N_\mathcal{F}^{\frac{p_J}{p_J'+\varepsilon}} \right)^\sharp (x) \leq \sup_{J: x \in J} \left( \frac{1}{|J|} \int_J N_{\mathcal{F}_J}(x)^{\frac{p_J}{p_J'+\varepsilon}} dx \right)$$

$$\leq b^{-p_J} \sup_{J: x \in J} \left( \frac{1}{|J|} \|Sf_J\|_{L^{p_J}(\mathbb{R}, l^{p_J'+\delta}(\mathcal{F}_J, L^1(J_T, l^2(T))))}^{p_J} \right)$$

$$\leq b^{-p_J} C\lambda^\varepsilon (1 + A^{-\frac{1}{\varepsilon}} \lambda)^{p_J} (M_{p_J}(Mf_J)(x))^{p_J}.$$

Again we can sharpen this argument to obtain

$$\left( N_\mathcal{F}^{\frac{p_J}{p_J'+\delta}} \right)^\sharp (x) \leq C b^{-p_J} \lambda^\varepsilon (1 + A^{-\frac{1}{\varepsilon}} \lambda)^{p_J} \max\{M_{p_J}(Mf_J)(x)^{p_J}, \lambda_0\}.$$

Taking the $\frac{p_J'+\delta}{p_J}$-norm on both sides proves estimate (60) and therefore also (42).

Georgia Institute of Technology, Atlanta, GA
*E-mail address*: lacey@math.gatech.edu

Christian-Albrechts-Universität at Kiel, Kiel, Germany
*E-mail address*: thiele@math.uni-kiel.de